\magnification=1200
\overfullrule=0pt
\centerline {\bf A structure theorem on non-homogeneous
linear equations in Hilbert spaces}\par
\bigskip
\bigskip
\centerline {BIAGIO RICCERI}\par
\bigskip
\bigskip
\noindent
{\bf Abstract:} A very
particular by-product of the result announced in the title reads as follows:
Let $(X,\langle\cdot,\cdot\rangle)$ be a real Hilbert space, $T:X\to X$ a
compact and symmetric
linear operator, and $z\in X$ such that the equation $T(x)-\|T\|x=z$ has
no solution
in $X$. For each $r>0$, set $\gamma(r)=\sup_{x\in S_r}J(x)$, where
$J(x)=\langle T(x)-2z,x\rangle$
and $S_r=\{x\in X:\|x\|^2=r\}$. Then,
the function $\gamma$ is $C^1$, increasing and strictly concave
in $]0,+\infty[$,
with $\gamma'(]0,+\infty[)=]\|T\|,+\infty[$; moreover, 
for each $r>0$, the problem of maximizing $J$ over
$S_r$ is well-posed,
and one has
$$T(\hat x_r)-\gamma'(r)\hat x_r=z$$
where $\hat x_r$ is the only global maximum of $J_{|S_r}$.\par
\par
\bigskip
\bigskip
\bigskip
Here and in the sequel, $(X,\langle\cdot,\cdot\rangle)$ is real Hilbert
space. For each $r>0$, set
$$S_r=\{x\in X : \|x\|^2=r\}\ .$$
In [1], we established the following result (with the usual
conventions $\sup\emptyset=-\infty$, $\inf\emptyset=+\infty$):\par
\medskip
THEOREM A ([1], Theorem 1). -  {\it  Let
$J:X\to {\bf R}$ be a sequentially weakly continuous $C^1$
functional, with $J(0)=0$.\par
Set
$$\rho=\limsup_{\|x\|\to +\infty}{{J(x)}\over {\|x\|^2}}$$
and
$$\sigma=\sup_{x\in X\setminus \{0\}}{{J(x)}\over {\|x\|^2}}\ .$$
Let $a, b$ satisfy
$$\max\{0,\rho\}\leq a<b\leq \sigma\ .$$
Assume that, for each $\lambda\in ]a,b[$, the
functional 
$x\to \lambda\|x\|^2-J(x)$
has a unique
global minimum, say $\hat y_{\lambda}$.
Let $M_a$ (resp. $M_b$ if $b<+\infty$ or
$M_b=\emptyset$ if $b=+\infty$) be
the set of all global minima of the functional
$x\to a\|x\|^2-J(x)$ (resp. $x\to b\|x\|^2-J(x)$
if $b<+\infty$).
Set
$$\alpha=\max\left \{0,\sup_{x\in M_b}\|x\|^2\right \}\ , $$
$$\beta=\inf_{x\in M_a}\|x\|^2 $$ 
and, for each $r>0$,
$$\gamma(r)=\sup_{x\in S_r}J(x)$$
Finally, assume that $J$ has no local maximum with norm less than $\beta$.\par
Then, the following assertions hold:\par
\noindent
$(a_1)$\hskip 5pt the function $\lambda\to g(\lambda):=
\|\hat y_{\lambda}\|^2$ is
decreasing in $]a,b[$ and its range is
$]\alpha,\beta[$\ ;\par
\noindent
$(a_2)$\hskip 5pt for each $r\in ]\alpha, \beta[$, the
point $\hat x_r:=\hat y_{g^{-1}(r)}$
is the unique global maximum of $J_{|S_r}$
and every maximizing sequence for
$J_{|S_r}$ converges to $\hat x_r$\ ;\par
\noindent
$(a_3)$\hskip 5pt 
the function $r\to \hat x_r$ is continuous in $]\alpha,\beta[$\ ;\par
\noindent
$(a_4)$\hskip 5pt
the function $\gamma$ is $C^1$, increasing
 and strictly concave in
$]\alpha,\beta[$\ ; \par
\noindent
$(a_5)$\hskip 5pt
 one has
$$J'(\hat x_r)=2\gamma'(r)\hat x_r$$
for all $r\in ]\alpha,\beta[$\ ;\par
\noindent
$(a_6)$\hskip 5pt one has
$$\gamma'(r)=g^{-1}(r)$$
for all $r\in ]\alpha,\beta[$.}\par
\medskip
We want to remark that, in the original statement of [1], one assumes that
$X$ is infinite-dimensional and that $J$ has no local maxima in $X\setminus
\{0\}$. These assumptions come from [2] whose
results are applied to get $(a_3)$, $(a_4)$ and $(a_5)$.
The validity of the current formulation just comes from the proofs themselves
given in [2] (see also [3]).\par
\smallskip
The aim of this very short paper is to show the impact of
Theorem A in the theory of non-homogeneous linear equations in $X$.\par
\smallskip
So, throughout the sequel,  $z$ is a non-zero
point of $X$ and $T:X\to X$ is a 
continuous linear operator.\par
\smallskip
We are interested in the study of the equation
$$T(x)-\lambda x=z$$
for $\lambda>\|T\|$. In this case, by the contraction mapping theorem, the
equation
has a unique non-zero solution, say $\hat v_{\lambda}$. Our structure
result just concerns such solutions.\par
\smallskip
As usual, we say that:\par
\noindent
-\hskip 5pt $T$ is compact if, for each bounded set $A\subset X$, 
 the set $\overline {T(A)}$ is compact\ ;\par
\noindent
-\hskip 5pt $T$ is symmetric if
$$\langle T(x),u\rangle=\langle T(u),x\rangle$$
for all $x,u\in X$\ .\par
\smallskip
We also denote by $V$ the set (possibly empty) of all solutions
of the equation
$$T(x)-\|T\|x=z$$
and set
$$\theta=\inf_{x\in V}\|x\|^2\ .$$
Of course, $\theta>0$.\par
\smallskip
Our result reads as follows:
\medskip
THEOREM 1. - {\it Assume 
that $T$ is compact and symmetric .\par
For each $\lambda>\|T\|$ and $r>0$, set
$$g(\lambda)=\|\hat v_{\lambda}\|^2$$
and
$$\gamma(r)=\sup_{x\in S_r}J(x)$$
where
$$J(x)=\langle T(x)-2z,x\rangle\ .$$
Then, the following assertions hold:\par
\noindent
$(b_1)$\hskip 5pt the function $g$ is decreasing in $]\|T\|,+\infty[$
and 
$$g(]\|T\|,+\infty[)=]0,\theta[\ ;$$
\noindent
$(b_2)$\hskip 5pt for each $r\in ]0,\theta[$, the point
$\hat x_r:=\hat v_{g^{-1}(r)}$ is the unique global
maximum of $J_{|S_r}$
and every maximizing sequence for
$J_{|S_r}$ converges to $\hat x_r$\ ;\par
\noindent
$(b_3)$\hskip 5pt the function $r\to \hat x_r$ is continuous
in $]0,\theta[$\ ;\par
\noindent
$(b_4)$\hskip 5pt the function $\gamma$ is $C^1$, increasing
and strictly concave in $]0,\theta[$\ ;\par
\noindent
$(b_5)$\hskip 5pt one has
$$T(\hat x_r)-\gamma'(r)\hat x_r=z$$
for all $r\in ]0,\theta[$\ ;\par
\noindent
$(b_6)$\hskip 5pt one has
$$\gamma'(r)=g^{-1}(r)$$
for all $r\in ]0,\theta[$\ .}\par
\medskip
Before giving the proof of Theorem 1, we establish the following\par
\medskip
PROPOSITION 1. - {\it Let $T$ be symmetric and let $J$ be defined as
in Theorem 1.\par
Then,
for $\tilde x\in X$, the following are equivalent:\par
\noindent
$(i)$\hskip 5pt $\tilde x$ is a local maximum of $J$\ .\par
\noindent
$(ii)$\hskip 5pt $\tilde x$ is a global maximum of $J$\ .\par
\noindent
$(iii)$\hskip 5pt $T(\tilde x)=z$ and $\sup_{x\in X}\langle T(x),x\rangle
\leq 0$\ .}\par
\smallskip
PROOF. First, observe that, since $T$ is symmetric,  the
functional $J$ is G\^ateaux differentiable and its derivative,
$J'$, is given by
$$J'(x)=2(T(x)-z)$$
for all $x\in X$ ([4], p. 235). By
 the symmetry of $T$ again, it is easy to check that,
for each $x\in X$, the inequality
$$J(\tilde x+x)\leq J(\tilde x) \eqno{(1)}$$
is equivalent to
$$\langle 2(T(\tilde x)-z)+T(x),x\rangle\leq 0\ .\eqno{(2)}$$
Now, if $(i)$ holds, then $J'(\tilde x)=0$ (that is $T(\tilde x)=z$) and
there is $\rho>0$ such that $(1)$ holds
for all $x\in X$ with $\|x\|\leq \rho$. So, from $(2)$, we have
$\langle T(x),x\rangle\leq 0$ for the same $x$ and then, by linearity,
 for all $x\in X$, getting $(iii)$. Vice versa, if $(iii)$ holds, then
 $(2)$ is satisfied for all $x\in X$ and so, by $(1)$, $\tilde x$ is a
global maximum of $J$, and the proof is complete.\hfill $\bigtriangleup$\par
\medskip
{\it Proof of Theorem 1}. For each $x\in X$, we clearly have
$$J(x)\leq \|T(x)-2z\|\|x\|\leq \|T\|\|x\|^2+2\|z\|\|x\|$$
and so
$$\limsup_{\|x\|\to +\infty}{{J(x)}\over {\|x\|^2}}\leq \|T\|\ .\eqno{(3)}$$
Moreover, if $v\in X\setminus \{0\}$ and $\mu\in {\bf R}\setminus
\{0\}$, we have
$${{J(\mu v)}\over {\|\mu v\|^2}}\geq -2{{\langle z,v\rangle}\over
{\mu\|v\|^2}}-\|T\|$$
and so
$$\limsup_{x\to 0}{{J(x)}\over {\|x\|^2}}=+\infty\ .\eqno{(4)}$$
Moreover, the compactness of $T$ implies that $J$ is sequentially weakly
continuous ([4], Corollary 41.9). Now, let $\lambda\geq\|T\|$. For each $x\in X$, set
$$\Phi(x)=\|x\|^2\ .$$
Then, for each $x,v\in X$, we have
$$\langle\lambda \Phi'(x)-J'(x)-(\lambda \Phi'(v)-J'(v)),x-v\rangle=
\langle 2\lambda (x-v)-2(T(x)-T(v)),x-v\rangle\geq $$
$$2\lambda\|x-v\|^2-2\|T(x)-T(v)\|\|x-v\|\geq 2(\lambda -\|T\|)\|x-v\|^2\ .\eqno{(5)}$$
 From $(5)$ we infer that the derivative of the functional $\lambda\Phi-J$
is monotone, and so the functional is convex. As a consequence, the critical
points of $\lambda\Phi-J$ are exactly its global minima. So, $\hat v_{\lambda}$ is
the only global minimum of $\lambda\Phi-J$ if $\lambda>\|T\|$ and $V$ is the set
of all global minima of $\|T\|\Phi-J$. Now, assume that $J$ has a
local maximum, say $w$. Then, by Proposition 1, $w$ is a global minimum of
$-J$ and $\sup_{x\in X}\langle T(x),x\rangle\leq 0$. Since $T$
is symmetric, this implies, in particular,
that $\|T\|$ is not in the spectrum of $T$. So, $V$ is a singleton.
By Proposition 1 of [1],
we have
$$\|w\|^2\geq \theta\ .$$
In other words, $J$ has no local maximum with norm less than $\theta$.
At this point, taking
$(3)$ and $(4)$ into account, we see that the assumptions
of Theorem A are satisfied (with $a=\|T\|$ and $b=+\infty$, and so
$\alpha=0$ and $\beta=\theta$),
and the conclusion follows directly from that result.\hfill $\bigtriangleup$
\medskip
Some remarks on Theorem 1 are now in order.\par
\medskip
REMARK 1. - Each of the two properties assumed on $T$ cannot be dropped.
Indeed, consider the following two counter-examples.\par
 Take $X={\bf R}^2$,
$z=(1,0)$ and $T(t,s)=(t+s,s-t)$ for all $(t,s)\in {\bf R}^2$.
So, $T$ is compact but not symmetric. In this case, we
have 
$$\hat x_r=(-\sqrt{r},0)\ ,$$
$$\gamma(r)=r+2\sqrt{r}$$
for all $r>0$. Hence, in particular, we have
$$T(\hat x_r)-\gamma'(r)\hat x_r=(1,\sqrt{r})\neq z\ .$$
That is, $(b_5)$ is not satisfied. \par
Now, take $X=l_2$, $z=\{w_n\}$, where $w_2=1$ and $w_n=0$ for
all $n\neq 2$,
 and $T(\{x_n\})=\{v_n\}$ for all
$\{x_n\}\in l_2$, where $v_1=0$ and $v_n=x_n$ for all $n\geq 2$.\par
So, $T$ is symmetric but not compact.
In this case, we have $\theta=+\infty$ and
$$\gamma(r)=r-2\sqrt{r}$$
for all $r\geq 4$. Hence, $\gamma$ is not strictly concave in $]0,+\infty[$.\par
\medskip
REMARK 2. - Note that the compactness of $T$ serves only to guarantee
that the functional $x\to \langle T(x),x\rangle$ is sequentially weakly
continuous. So, Theorem 1 actually holds under such a weaker condition.
\medskip
REMARK 3. - A natural question is: if 
assertions $(b_1)-(b_6)$ hold, must the operator $T$ be
symmetric
and the functional $x\to \langle T(x),x\rangle$ sequentially
weakly continuous ?\par
\medskip
REMARK 4. - Note that if $T$, besides to be compact
and symmetric, is also positive (i.e. $\inf_{x\in X}\langle T(x),x\rangle\geq
0$), 
then, by classical results, the operator $x\to T(x)-\|T\|x$ is not
surjective, and so there
are $z\in X$ for which the conclusion of Theorem 1 holds
with $\theta=+\infty$.\par
\medskip
We conclude with an application Theorem 1 to a classical
Dirichlet problem.\par
\smallskip
So, let $\Omega\subset {\bf R}^n$ be a bounded domain with
smooth boundary. Let $\lambda_1$ be the first eigenvalue
of the problem
$$\cases {-\Delta u=\lambda u & in $\Omega$\cr & \cr
u=0 & on $\partial \Omega$\ .\cr}$$
Fix a non-zero continuous function $\varphi:\overline {\Omega}\to
{\bf R}$.\par
\smallskip
For each $\mu\in ]0,\lambda_1[$, let $u_{\mu}$ be the
unique classical solution of the problem
$$\cases {-\Delta u=\mu(u+\varphi(x)) & in $\Omega$\cr & \cr
u=0 & on $\partial \Omega$\ .\cr}$$
Also, set
$$\psi(\mu)=\int_{\Omega}|\nabla u_{\mu}(x)|^2dx$$
and 
$$\eta(r)=\sup_{u\in U_r}\Phi(u)$$
where
$$\Phi(u)=\int_{\Omega}|u(x)|^2dx+2\int_{\Omega}\varphi(x)u(x)dx$$
and
$$U_r=\left \{ u\in H^1_0(\Omega) : \int_{\Omega}|\nabla u(x)|^2dx=r\right \}\ .$$
Finally, denote by $A$ the set of all classical solutions of the
problem
$$\cases {-\Delta u=\lambda_1(u+\varphi(x)) & in $\Omega$\cr & \cr
u=0 & on $\partial \Omega$\cr}$$
and set
$$\delta=
\inf_{u\in A}\int_{\Omega}|\nabla u(x)|^2dx \ .$$
Then, by using standard variational methods, we can directly draw the following
result from Theorem 1 :\par
\medskip
THEOREM 2. - {\it The following assertions hold:\par
\noindent
$(c_1)$\hskip 5pt the function $\psi$ is increasing in $]0,\lambda_1[$ and
one has
$$\psi(]0,\lambda_1[)=]0,\delta[\ ;$$
\noindent
$(c_2)$\hskip 5pt for each $r\in ]0,\delta[$, the function $w_r:=u_{\psi^{-1}(r)}$ is the unique
global maximum of $\Phi_{|U_r}$ and each maximizing sequence for
$\Phi_{|U_r}$ converges to $w_r$ with respect to the topology of $H^1_0(\Omega)$\ ;\par
\noindent
$(c_3)$\hskip 5pt the function $r\to w_r$ is continuous in $]0,\delta[$ with respect to the
topology of $H^1_0(\Omega)$\ ;\par
\noindent
$(c_4)$\hskip 5pt the function $\eta$ is $C^1$, increasing and strictly concave in $]0,\delta[$\ ;\par
\noindent
$(c_5)$\hskip 5pt for each $r\in ]0,\delta[$, the function $w_r$ is the unique classical solution of the problem
$$\cases {-\Delta u={{1}\over {\eta'(r)}}(u+\varphi(x)) & in $\Omega$\cr & \cr
u=0 & on $\partial \Omega$\ ;\cr}$$
\noindent
$(c_6)$\hskip 5pt one has
$$\eta'(r)={{1}\over {\psi^{-1}(r)}}$$
for all $r\in ]0,\delta[$.}\par
\vfill\eject
\centerline {\bf References}\par
\bigskip
\bigskip
\noindent
[1]\hskip 5pt B. RICCERI, {\it On a theory by Schechter and Tintarev}, Taiwanese J. Math., {\bf 12}
(2008), 1303-1312.
\smallskip
\noindent
[2]\hskip 5pt M. SCHECHTER and K. TINTAREV, {\it Spherical maxima in
Hilbert space and semilinear elliptic eigenvalue problems}, Differential
Integral Equations, {\bf 3} (1990), 889-899.\par
\smallskip
\noindent
[3]\hskip 5pt K. TINTAREV, {\it Level set maxima and quasilinear elliptic problems},
Pacific J. Math., {\bf 153} (1992), 185-200.\par
\smallskip
\noindent
[4]\hskip 5pt E. ZEIDLER, {\it Nonlinear functional analysis and its
applications}, vol. III, Springer-Verlag, 1985.\par
\bigskip
\bigskip
\bigskip
\bigskip
Department of Mathematics\par
University of Catania\par
Viale A. Doria 6\par
95125 Catania\par
Italy\par
{\it e-mail address}: ricceri@dmi.unict.it

\bye